\documentclass[12pt]{article}

\usepackage[reqno,tbtags]{amsmath}

\usepackage{amssymb}

\def\thebibliography#1{\vspace{0.5cm} {\flushleft \bf
References \\ } \list
 {[\arabic{enumi}]}{\settowidth\labelwidth{[#1]}\leftmargin\labelwidth
 \advance\leftmargin\labelsep
 \usecounter{enumi}}
 \def\newblock{\hskip .11em plus .33em minus -.07em}
 \sloppy
 \sfcode`\.=1000\relax}

\protect\newcounter{theoremnumber}

\protect\newcounter{corollarynumber}

\protect\newcounter{remarknumber}

\protect

\newskip\halflineskip

\newenvironment{theorem}%
{\par\addvspace{\halflineskip} \refstepcounter{theoremnumber}
{\noindent \bf Theorem~\thetheoremnumber .} %
\hskip.5em\ignorespaces\it}%
{\vskip\halflineskip\par\rm}

\newenvironment{corollary}%
{\par\addvspace{\halflineskip} \refstepcounter{corollarynumber}
{\noindent \bf Corollary~\thecorollarynumber .}%
\hskip.5em\ignorespaces\it}%
{\vskip\halflineskip\par\rm}

{\par\addvspace{\halflineskip} \refstepcounter{theoremnumber}
{\noindent \bf Proposition~\thetheoremnumber .}%
\hskip.5em\ignorespaces\it}%
{\vskip\halflineskip\par\rm}

\newenvironment{lemma}%
{\par\addvspace{\halflineskip} \refstepcounter{theoremnumber}
{\noindent \bf Lemma~\thetheoremnumber .}%
\hskip.5em\ignorespaces\it}%
{\vskip\halflineskip\par\rm}

{\par\addvspace{\halflineskip} \refstepcounter{definitionnumber}
{\noindent \bf Definition~\thedefinitionnumber .}%
\hskip.5em\ignorespaces}%
{\vskip\halflineskip\par}

\newenvironment{example}%
{\par\addvspace{\halflineskip} \refstepcounter{sergei}
{\noindent \bf Example~\thesergei .}%
\hskip.5em\ignorespaces}%
{\vskip\halflineskip\par}

\newenvironment{remark}%
{\par\addvspace{\halflineskip} \refstepcounter{remarknumber}
{\noindent \bf Remark~\theremarknumber .}%
\hskip.5em\ignorespaces}%
{\vskip\halflineskip\par}

{\par\addvspace{\halflineskip}}%
{\vskip\halflineskip\par}

\def\@begintheorem#1#2{\par\addvspace{\halflineskip}
{\bf #1\  #2.} \hskip.5em}

\def\@opargbegintheorem#1#2#3{\par\vskip\halflineskip
 {\bf #1\ #2.\ {\rm (#3)}.}\hskip.5em}

\def\@endtheorem{\vskip\halflineskip\par}

\newenvironment{proof}{\noindent{\bf Proof. }\rm}
{\unskip\nobreak\hfil\penalty50\hskip1em\hbox{}
\nobreak\hfill\qed\par\smallskip}

\def\qed{\vrule height1ex width1ex depth0pt}

{{\raisebox{0.4ex}{$\bigtriangledown$}}\vskip\halflineskip\par}

\protect

{\protect\nopagebreak\section*{\protect\raggedright \small \it
Acknowledgements}}%
{}

{\protect\nopagebreak\section*{\protect\raggedright \normalsize
Acknowledgement}}%
{}

\protect

\protect\newcommand{\sectionhead}[1]%
{\section {#1}}

\def\abstract{\vspace{0cm}
{\bf \small Abstract.}\footnotesize }

\begin{document}

\title{ Monotone operator functions on $C^*$-algebra}

\author{Hiroyuki Osaka \vspace{0.3cm}\\
Department of Mathematical Sciences, Ritsumeikan University,\\ 
Kusatsu, Shiga 525-8577, Japan \\
e-mail: osaka@se.ritsumei.ac.jp \\
FAX:   +81 77 561 2657 tel: +81 77 561 2656 \vspace{0.5cm}
\\ Sergei D. Silvestrov \vspace{0.3cm}\\
Centre for Mathematical Sciences, Department of Mathematics, \\
Lund Institute of Technology,
Box 118, SE-22100 Lund, Sweden. \\
e-mail: sergei.silvestrov@math.lth.se \\
FAX:  +46 46 2224010 \ \ tel: +46 46 2228854 \vspace{0.5cm}
\\ 
Jun Tomiyama \vspace{0.3cm}\\ 
Prof.Emeritus of Tokyo Metropolitan University,\\
201 11-10 Nakane 1-chome, \\ 
Meguro-ku, Tokyo, Japan \\
e-mail: jtomiyama@fc.jwu.ac.jp
}
\date{\today}
%\begin{document}
\maketitle

%\pagebreak

%

%

%\vspace*{7.5cm}

\begin{abstract}
The article is devoted to investigation of classes of functions
monotone as functions on general $C^*$-algebras that are not
necessarily the $C^*$-algebras of all bounded linear operators on
a Hilbert space as it is in classical case of matrix and operator
monotone functions. We show that for general $C^*$-algebras the
classes of monotone functions coincide with the standard classes
of matrix and operator monotone functions. For every class we give
exact characterization of $C^*$-algebras that have this class of
monotone functions, providing at the same time a monotonicity
characterization of subhomogeneous $C^*$-algebras. We use this
characterization to generalize one function based monotonicity
conditions for commutativity of a $C^*$-algebra, to one function
based monotonicity conditions for subhomogeneity. As a
$C^*$-algebraic counterpart of standard matrix and operator
monotone scaling, we investigate, by means of projective
$C^*$-algebras and relation lifting, the existence of
$C^*$-subalgebras of a given monotonicity class.

\end{abstract}

\footnotetext{Mathematics Subject Classification 2000: 46L05
% Primary
% 47L65; Secondary 37B05, 37B20, 47L55
} \footnotetext{This work was supported by The Royal Swedish
Academy of Sciences and by Crafoord foundation and 
JSPS Grant for Scientific Research No. 14540217(c)(1).}

\section{Introduction.} \label{sec:intr}

The real-valued continuous function $f: I \mapsto \mathbb{R}$ on a
(non trivial) interval $I \neq \mathbb{R}$ is called $A$-monotone
for a given $C^*$-algebra $A$ if for any $x, y \in A$ with
spectrum in $I$, \

\begin{equation} \label{ineq:monot}
x\leq_A y \quad \Rightarrow \quad f(x)  \leq f(y)
\end{equation}

We denote by $P_A(I)$ the set of all $A$-monotone functions
(defined on the interval I) for a $C^*$-algebra $A$. 
If $A = B(H)$, 
the standard $C^*$-algebra of all bounded linear operators on a
Hilbert space $H$, then $P_A(I)=P_{B(H)}(I)$ is called the set of all
operator monotone functions. If $A = M_n $, the standard
$C^*$-algebra of all complex $n\times n$ matrices or equivalently
of all (bounded) linear operators on an $n$-dimensional complex
Hilbert space, then $P_{n}(I)=P_A(I)=P_{M_n}(I)$ is called the set of
all matrix monotone functions of order $n$ on an interval $I$. The
set $P_{n}(I)$ consists of continuous functions on $I$ satisfying
\eqref{ineq:monot} for pairs $(x,y)$ of self-adjoint $n\times n$
matrices with spectrum in $I$. 
For each positive integer $n$, the
proper inclusion $P_{n+1}(I) \subsetneq P_{n}(I)$ holds
\cite{Donoghuebook,HansenJiTomiyama-art}. 
For infinite-dimensional
Hilbert space, the set of operator monotone functions on $I$ can
be shown to coincide with the intersection
$$P_{\infty}(I) = \bigcap_{n=1}^{\infty} P_n(I),$$ or in other
words a function is operator monotone if and only if it is matrix
monotone of order $n$ for all positive integers $n$ 
\cite[Chap.5, Prop.5.1.5 (1)]{HiaiYanagibook}. 
Keeping this in mind, for
infinite-dimensional Hilbert space, we denote the class of
operator monotone functions also by $P_{\infty}(I)$ or simply by
$P_{\infty}$ when the choice of the interval is clear from
context. 
For the sake of clarity, if not stated otherwise we will
assume that all $C^*$-algebras contain a unity. Formulations of
most of the results can be adjusted to hold also in non-unital
situation by the standard procedure of adjoining the unity, which
amounts to adding a one-dimensional irreducible representation to
the set of irreducible representations of $A$.

The Section \ref{sec:scalingtheorems} is devoted to description of
the classes of monotone operator functions of $C^*$-algebras. We
show that for general $C^*$-algebras the classes of monotone
functions are the standard classes of matrix and operator monotone
functions. For every such class we give exact characterization of
$C^*$-algebras that have this class of monotone functions. This
can be also used to give a monotonicity characterization of
subhomogeneous $C^*$-algebras as discussed in \cite[Theorem 5]{HansenJiTomiyama-art}. In Section
\ref{sec:moncharcomsubhom} we use these characterizations to
generalize one function based monotonicity condition for
commutativity of a $C^*$-algebra, obtained by T.~Ogasawara
\cite{Ogasawara} and G.~K.~Pedersen \cite{Pedersenbook}, W. Wu
\cite{Wu-art}, and Ji and Tomiyama \cite{JiTomiyama}, to one
function based monotonicity condition for subhomogeneity. Finally
in Section \ref{sec:existsubalg}, we investigate, as a
\mbox{$C^*$-algebraic} counterpart of standard matrix and operator
monotone scaling, the existence of $C^*$-subalgebras of a given
monotonicity class. We also state several problems motivated by
the obtained results.

\section{Scaling theorems} \label{sec:scalingtheorems}

To begin with, note that for any $C^*$-algebra $A$ there is a
Hilbert space $H$ such that $P_{B(H)} \subseteq P_A$, and in
particular always $P_\infty \subseteq P_A$. Indeed, by
Gelfand-Naimark construction $A$ is isometrically isomorphic to a
$C^*$-subalgebra $\tilde{A}$ of $B(H)$ for some Hilbert space $H$.

Any isomorphism between two $C^*$-algebras preserves the standard
partial order induced by their positive cones. 
Therefore, any function
which is operator monotone, that is $B(H)$-monotone, is also
$\tilde{A}$-monotone and hence $A$-monotone.  In general, if
$B\hookrightarrow A$ that is a $C^*$-algebra $B$ is isomorphic to
a $C^*$-subalgebra of a $C^*$-algebra $A$, 
then $P_{A} \subseteq P_B$. 
In other words the mapping $A \mapsto P_{A}$ is
non-increasing. For the standard matrix imbedding scaling we have

$$
M_1 \hookrightarrow M_2 \hookrightarrow M_3
\hookrightarrow \dots \hookrightarrow M_k \hookrightarrow \dots
\hookrightarrow B(H).$$ 
This standard imbedding  sequence is
infinite and strictly increasing if $\dim  H = \infty$, and we
have the corresponding decreasing sequence 
$$
P_{1}(I) \supset
P_{2}(I) \supset P_{3}(I) \supset \dots \supset P_{n}(I) \supset
\dots \supset P_\infty (I).
$$ 
The inclusions of function spaces
$P_{n+1}(I) \subset P_{n}(I)$ and $P_{\infty}(I) \subset P_{n}(I)$
are strict for all positive integers $n$ and non-trivial intervals $I$. 
Even though this fact has been known almost from the
beginning of the theory of operator monotone functions, only
recently explicit examples of functions from $P_{n} \setminus
P_{n+1}$ for arbitrary choice of $n$ have been constructed
\cite{HansenJiTomiyama-art}. For general $C^*$-algebras, the
imbedding partial order is more flexible allowing for different
kinds of scalings.

The irreducible representations contain an important information
about $C^*$-algebras, and dimensions of representations are the 
important classifying parameter. 
The following Lemma, which sharpens the assertion of 
\cite[Theorem 5,(1) and (2)]{HansenJiTomiyama-art},
is a key to
further understanding of relationship between dimensions of
irreducible representations and $A$-monotonicity for a
$C^*$-algebra $A$ on one side, and the operator monotonicity and
matrix monotonicity on the other.

In the sequel, without loss of generality
\cite{HansenJiTomiyama-art}, we assume that $I = [0,\infty[$ and
drop the interval from the corresponding notations.

\begin{lemma} \label{th:mainlemmareps-omf}
Let $A$ be a (unital) $C^*$-algebra.
 \begin{itemize}
 \item[1)] If $A$ has an irreducible representation of dimension $n$
 then any $A$-monotone function becomes
$n$-matrix monotone, that is $P_A \subseteq P_n $.
 \item[2)] If $\dim \pi \leq n$ for any irreducible representation
 $\pi$ of $A$, then $P_n \subseteq P_A$.
 \item[3)] If the set of dimensions of finite-dimensional irreducible representa tions of $A$ is unbounded, then every $A$-monotone function
is operator monotone, that is $P_A = P_\infty$.
 \item[4)] If $A$ has an infinite-dimensional irreducible
 representation, then every $A$-monotone function
 is operator monotone, that is $P_A = P_\infty$.
\end{itemize}

\end{lemma}

\begin{proof}
Let $\pi:A\rightarrow M_n$ be an $n$-dimensional irreducible
representation of $A$. Then irreducibility implies that $\pi (A) =
M_n$. Thus for any pair $c,d \in M_n$, such that $0\leq c \leq d$
there exists $a, b \in A$ such that $0\leq a \leq b$ and 
$\pi(a) = c$ and $\pi(b) = d$. Then $f(a) \leq f(b)$ and hence 
$\pi(f(a)) \leq \pi(f(b))$ for any $f\in P_A(I)$. By continuity, 
$\pi (f(x)) = f( \pi (x))$ for any $x\in A$. 
Thus $f(c) = f( \pi (a)) \leq f(\pi (b)) = f(d)$, 
and therefore $f\in P_n$. Hence, we have proved that
$P_A \subseteq P_n$.

2) For any $f \in P_n$, for any $0\leq a \leq b$ in $A$ and for
any irreducible representation $\pi: A \rightarrow M_m$, where 
$m \leq n$, we have $\pi (a) \leq \pi (b)$ in $M_m$. 
Then $\pi(f(a)) = f(\pi (a)) \leq f (\pi (b)) = \pi(f(b))$. 
If $0\leq \pi(f(b)-f(a))$ for any irreducible representation $\pi$, then
$spec(f(b)-f(a))\in [0,\infty [ $ that is $0\leq f(b)-f(a) $ or
equivalently $f(a)\leq f(b)$. Thus, $f\in P_A$ and we proved that
$P_n \subseteq P_A$.

3) Let $\{\pi_j \mid j\in \mathbb{N}\setminus \{0\}\}$ be a
sequence of irreducible finite-dimensional representations of $A$
such that $n_j= \dim \pi_j \rightarrow \infty$ when $j\rightarrow
\infty$. By 1) we have inclusion $P_A \subseteq P_{n_k}$ for any
$k\in \mathbb{N}\setminus \{0\}$. Hence 
$$
P_A \subseteq \bigcap_{k\in \mathbb{N}\setminus \{0\}} P_{n_k} 
= \bigcap_{k\in\mathbb{N}\setminus \{0\}} P_{k} = P_\infty,
$$ and since always
$P_\infty \subseteq P_A$ holds, we get the equality 
$P_A = P_\infty$.

4) Let $\pi:A\rightarrow B(H)$ be irreducible representation of
$A$ on an infinite-dimensional Hilbert space $H$. By Kadison
transitivity theorem, in the form it is stated in Takesaki's book
\cite[Ch.2, Theorem 4.18]{takesaki-bok}, $\pi (A) p = B(H) p$ for
every projection $p: H\rightarrow H $ of a finite rank 
$n=\dim p H < \infty $. 
Let $B=\{a\in A \mid \pi (a) pH \subseteq pH, \pi(a)^*
pH \subseteq pH \}$ be the $C^*$-subalgebra of $A$ consisting of
elements mapped by $\pi$ to operators that, together with their
adjoints, leave $pH$ invariant. The restriction of 
$\pi: B \mapsto p B(H) p $ to $B$ is $n$-dimensional representation of $B$ on
$pH$, and moreover it is is irreducible and surjection, 
since $\pi(B) p = p \pi (B) p = p \pi (A) p = p B(H) p = B (pH)$. Thus 1)
yields $P_A \subseteq P_B \subseteq P_n$, since $B$ is a
$C^*$-subalgebra of $A$. As the positive integer $n$ can be chosen
arbitrary, we get the inclusion 
$$
P_A \subseteq \bigcap_{n\in\mathbb{N} \setminus \{0\}} P_n = P_\infty .
$$ 
Combining it with
$P_\infty \subseteq P_A$ yields the equality $P_A = P_\infty$.
\end{proof}

\begin{corollary} \label{th:maxmonmaxhomog}
If $n_0 = \sup \{k \mid P_A \subseteq P_{k}\}$, then
$$
n_0 = n_1 =\sup \{\dim (\pi) \mid \pi \text{ is irreducible representation
of } A \}.
$$
\end{corollary}

\begin{proof}
 By Lemma \ref{th:mainlemmareps-omf}, the positive integer
 $n_0 = \sup \{k  \mid P_A \in P_{k}\}$
 exists only if the set of dimensions of irreducible
 representations of $A$ is bounded. Let
 $$
 n_1 = \sup \{\dim (\pi) \mid \pi \text{ is irreducible representation
of } A \}.
$$ 
Then by 1) and 2) of Lemma \ref{th:mainlemmareps-omf}
we have $P_{n_1} \subseteq P_A \subseteq P_{n_1}$, and hence 
$P_A = P_{n_1}$. Thus $P_{n_1} = P_A \subseteq P_{n_0}$ by 1) of Lemma
\ref{th:mainlemmareps-omf}. So, $n_1 \geq n_0$, and since 
$n_0 = \sup \{n \mid P_A \subseteq P_n \}$ and $P_A= P_{n_1}$, we get
the desired $n_0=n_1$. If $n_0 = \infty$, then $P_A = P_\infty$.
By Lemma \ref{th:mainlemmareps-omf} either $A$ has an
infinite-dimensional irreducible representation or the set of
dimensions of irreducible representations is unbounded, that is
$n_0=n_1 = \infty$, because if on the contrary the set of
dimensions of irreducible representations is bounded by some
positive integer $n$, then $P_n \subseteq P_A = P_\infty$, which
is impossible since $P_\infty \subset P_n$ with gap $P_n \setminus
P_\infty \neq \emptyset$.
\end{proof}

Recall that a $C^*$-algebra $A$ is said to be subhomogeneous if the set of 
dimensions of its irreducible representations is bounded. We say 
that $A$ is $n$-subhomogeneous or subhomogeneous of degree $n$ if 
$n$ is the highest  
dimension of those irreducible representations of $A$.

%\pagebreak

\begin{theorem}\label{th:monot-homog} 
Let $A$ be a $C^*$-algebra. Then

 \begin{itemize}
  \item[1)] $P_A = P_{\infty}$ if and only if
either the set of dimensions of finite-dimensional irreducible representations 
of $A$ is unbounded,
or $A$ has an infinite-dimensional irreducible representation.

  \item[2)] $P_A = P_n$ for some positive integer $n$ if and only if $A$ is
$n$-subhomogeneous.
\end{itemize}
\end{theorem}

\begin{proof}
By Lemma \ref{th:mainlemmareps-omf} the only part of 1) left to
prove is that $P_A = P_{\infty}$ implies that either the set of
dimensions of finite-dimensional irreducible representations of
$A$ is unbounded, or $A$ has an infinite-dimensional irreducible
representation. Suppose on the contrary that 
$$
n_1 = \sup \{\dim
(\pi) \mid \pi \text{ is irreducible representation of } A \} <
\infty.
$$ 
Then $P_A \subseteq P_{n_1}$ by Corollary
\ref{th:maxmonmaxhomog}, and $P_{n_1} \subseteq P_A$ by 2) of
Lemma \ref{th:mainlemmareps-omf}. Thus $P_A = P_{n_1} $. 
But there
is a gap between $P_\infty$ and $P_n$ for any $n$. 
Hence $P_A \neq P_\infty$, in contradiction to the initial assumption 
$P_A = P_{\infty}$.

In part 2), again thanks to Lemma \ref{th:mainlemmareps-omf}, it
is left to prove that if $P_A = P_n$, then $A$ is
$n$-subhomogeneous. If $P_A=P_n$, then 
$$
n=n_0 = \sup \{k\in
\mathbb{N} \mid P_A \subseteq P_{k}\}.
$$ 
Indeed, if $n_0 > n$, then 
$ P_{n_0} \subsetneq P_A =P_n $ since there exists a gap 
$P_m \subsetneq P_n$ for all $m > n$ as proved 
in \cite{HansenJiTomiyama-art}. 
But
this contradicts to $P_A \subseteq P_{n_0}$ true by definition of
$n_0$. 
Hence $n_0 \leq n$. By Corollary \ref{th:maxmonmaxhomog},
$$
n_0 =n_1 = \sup \{\dim (\pi) \mid \pi \text{ is irreducible
representation of } A \},$$ 
and thus $n\leq n_0$. Therefor, 
$n = n_0 = n_1$ and so $A$ is $n$-subhomogeneous.
\end{proof}

\begin{remark} \label{rem:PAPkinclud} A useful observation is that by Lemma
\ref{th:mainlemmareps-omf} and Theorem \ref{th:monot-homog}, for
any $C^*$-algebra $A$ and any positive integer $k$, only two
possibilities are possible, either $P_A \cap P_k = P_k$ or $P_A
\cap P_k = P_A$.
\end{remark}

The Theorem \ref{th:monot-homog} can be used to extend 2) of Lemma
\ref{th:mainlemmareps-omf} to be the "if and only if " statement
also proved in \cite{HansenJiTomiyama-art}.

\begin{corollary} \label{Cor:maxmonmaxhomog}
Every matrix monotone function of order $n$ is $A$-monotone if and
only if the dimension of every irreducible representation of $A$
is less or equal to $n$.
\end{corollary}

\begin{proof}

The "if" part is 2) of Lemma \ref{th:mainlemmareps-omf}. To prove
the "only if" part note that $P_n \subseteq P_A$ implies that 
$P_A = P_m$ for some $m\leq n$, and by 2) of Theorem
\ref{th:monot-homog}, $m = n_1$. 
Thus $\dim (\pi) \leq n_1 = m \leq n $ 
for any irreducible representation of $A$.
\end{proof}

If  $A$ is a commutative $C^*$-algebra, then every irreducible
representation of $A$ is one-dimensional and hence $P_A = P_1$,
the set of all non-decreasing continuous functions. A natural
class generalizing commutative $C^*$-algebras consists of
$n$-homogeneous $C^*$-algebras, that is $C^*$-algebras with all
non-zero irreducible representations being $n$-dimensional. The
description of the class of monotone operator functions for
$n$-homogeneous $C^*$-algebras follows from Theorem
\ref{th:monot-homog}.

\begin{corollary} \label{Cor:nhomogalgMOF}
If $C^*$-algebra $A$ is $n$-homogeneous, then $P_A = P_n$.
\end{corollary}

\begin{example}
For any Hilbert space $H$, the equality $P_{B(H)} = P_{\dim H}$
holds. If $\dim H =n < \infty$, then $B(H) = M_n$ and 
$P_{B(H)} = P_n$; and if $\dim H = \infty$, then $P_{B(H)} = P_\infty$.
\end{example}

\begin{example}
The irrational rotation $C^*$-algebra $A_\theta$ is the
$C^*$-algebra generated by two unitaries $u$ and $v$ satisfying
commutation relation $uv =e^{i 2 \pi \theta}vu$ with some
irrational $\theta \in ]0,1[ $. It is isomorphic to the crossed
product $C^*$-algebra $C(\mathbb{T}) \rtimes_{\sigma_\theta}
\mathbb{Z}$ associated to the dynamical system consisting of the
rotation $\sigma_\theta$ of the one-dimensional torus (the unite
circle) $\mathbb{T}$ by an angle $2 \pi \theta$ with irrational
$\theta$. All non-zero irreducible representations of $A_\theta$
are infinite-dimensional since all points of $\mathbb{T}$ are
aperiodic under action of $\sigma_\theta$. 
Hence $P_{A_\theta}= P_\infty $ by Theorem \ref{th:monot-homog}.

The rational rotation $C^*$-algebra $A_\theta$ is the crossed
product $C^*$-algebra $C(\mathbb{T}) \rtimes_{\sigma_\theta}
\mathbb{Z}$ where $\sigma_\theta$ is the rotation of $\mathbb{T}$
by the angle $2 \pi \theta$ with rational $\theta = \frac{m}{n}$
($m$ and $n$ are relatively prime). The $C^*$-algebra $A_\theta$
is isomorphic to the $C^*$-algebra of cross-sections in the fibre
bundle over $\mathbb{T}^2$ with fibre $M_n$, the $n\times n$
matrix algebra, and the structure group $U_n$, the $n$-dimensional
unitary group. All points of $\mathbb{T}$ are periodic of period
$n$ and thus all irreducible representations of $A_\theta$ are
$n$-dimensional, which means that $A_\theta$ is an $n$-homogeneous
$C^*$-algebra. Hence $P_{A_\theta}= P_n $ by Corollary
\ref{Cor:nhomogalgMOF}.
\end{example}

\begin{example}
The $C^*$-algebras $C_{0}(X, M_n (\mathbb{C}))$ of continuous and
vanishing at infinity $M_n (\mathbb{C})$-valued functions on a
locally compact Hausdorff space $X$ are $n$-homogeneous
$C^*$-algebras, and hence Hence
$P_{C_{0}(X, M_n (\mathbb{C}))}= P_n $ by Corollary \ref{Cor:nhomogalgMOF}. 
These $C^*$-algebras
can be viewed as the space of continuous sections, vanishing at
infinity, of the trivial $M_n$-bundle $X \times M_n$. In fact
every $n$-homogeneous algebra arises as the algebra of continuous
sections of some $M_n$-bundle \cite{Fellart1,TomTohokuJ62,TomTak61}.
\end{example}

\begin{example}
Let $A=C(\mathbb{T}) \rtimes_\sigma \mathbb{Z}$ be the crossed
product algebra associated to the dynamical system consisting of a
homeomorphism $\sigma$ of $\mathbb{T}$. If $\sigma$ is an
orientation preserving homeomorphism of the circle  without
periodic points, then $n_1 = \infty$ and hence $P_A = P_\infty$ by
Theorem \ref{th:monot-homog}.
\end{example}

\begin{example}
Let $A=C(X) \rtimes_\sigma \mathbb{Z}$ be the crossed product
algebra associated to the dynamical system consisting of a
homeomorphism $\sigma$ of a compact Hausdorff space $X$. Then
since any finite-dimensional irreducible representation of $A$ is unitarily
equivalent to an induced representation arising from a periodic
point according to \cite[Proposition 4.5]{TomSeoulLN1} (see also
\cite{SilTomTypeIart}), the equality $P_A= P_\infty$ holds if and
only if $(X, \sigma)$ either have an aperiodic orbit or the set of
periods of periodic points in $(X, \sigma)$ is an unbounded subset of positive
integers; and if all points of $X$ are periodic for $\sigma$, and
the set of periods is bounded, then $P_A = P_n$ for the maximal
period $n$, coinciding with the maximal dimension for irreducible
representations of $A$.
\end{example}

\begin{example} \label{ex:Heisgroup}
Let ${\cal H}$ be the three-dimensional discrete Heisenberg group
represented by matrices,
$${\cal H} =\left\{\left(\begin{array}{ccc} 1, l, m \\

                                     0, 1, n \\

                                     0, 0, 1

                  \end{array}\right) \mid l,m,n \in \mathbb{Z} \right\}.$$

Then it can be shown that the group $C^*$-algebra $C^*({\cal H})$
is isomorphic to the crossed product $C^*$-algebra 
$C^*({\cal H}) = C(\mathbb{T}^2) \rtimes_\sigma \mathbb{Z}$ 
associated to the
homeomorphism of the two-dimensional torus $\mathbb{T}^2$ defined
by $\sigma (s,t) = (s, t-s)$. This homeomorphism acts as rational
rotation along the second coordinate direction if $s$ is rational,
and as irrational rotation if $s$ is irrational. This means in
particular that $C^*({\cal H})$ has irreducible representations of
infinite dimension and of any finite dimension. Hence
$P_{C^*({\cal H})} = P_\infty$ by Theorem \ref{th:monot-homog}.
\end{example}

\begin{example}
The $c_0$-direct sum $A=\sum_{i=1}^\infty \oplus M_{n_i}$ of
matrix algebras with a sequence of dimensions such that 
$n_i \rightarrow \infty$ when $i\rightarrow \infty$ is an example of a
$C^*$-algebra for which there all irreducible representations are
finite-dimensional, but the set of dimensions is unbounded. For
this $C^*$-algebra $P_A = P_\infty$ by Theorem \ref{th:monot-homog}.
\end{example}

\begin{example}

For any positive integer $n\geq 2$, the Cuntz $C^*$-algebra 
${\cal O}_n$ is the universal unital $C^*$-algebra on generators
$s_1,\dots, s_n$ satisfying relations

\begin{eqnarray*} && s_1s^*_1 + \dots + s_n
s^*_n = 1 \\
&& s^*_js_k = \delta_{jk}1 = \left\{\begin{array}{l} 1, \mbox{ if
}
j = k \\

0, \mbox{ if } j \neq k \end{array} \right. \mbox{ for } j,k
= 1,\dots, n.
\end{eqnarray*}

As ${\cal O}_n$ has infinite-dimensional irreducible
representations, $P_{{\cal O}_n} = P_\infty$ by 4) of Lemma
\ref{th:mainlemmareps-omf}. 
Actually ${\cal O}_n$ is known to be an infinite-dimensional simple 
C*-algebra.
\end{example}

\section{Monotonicity characterizations of commutativity and sub-homogeneity}
\label{sec:moncharcomsubhom}

In $C^*$-algebras order induced by positivity is closely connected
to algebraic properties. As an outcome of this, one can prove
several unexpected results characterizing such properties as
commutativity and sub-homogeneity in terms of monotonicity
properties of functions.

There are several characterizations for the commutativity of
$C^*$-algebras.  One type  is the well-known  Stinespring theorem,
that is, a C$^*$-algebra $A$ is commutative if and only if  every
positive linear map from $A$ to another C$^*$-algebra $B$ (or from
$B$ to $A$) becomes completely positive \cite{Stinespring-art}. To
be precise, $A$ becomes commutative if and only if every positive
linear maps to $B$ becomes two-positive (and then automatically
completely positive). This is the beginning of the long and
fruitful developments of understanding the matricial order
structure of operator algebras (see for example
\cite{EffrosRuanbook}).

One of the first results in the direction of operator algebraic
monotonicity of functions, obtained in 1955 by T. Ogasawara
\cite{Ogasawara}, states that if $0\leq x \leq y$ implies $x^2
\leq y^2$ for all $x,y$ in a $C^*$-algebra $A$, then $A$ is
commutative. A proof of this result can be found also in the G. K.
Pedersen's book 
\cite[Proposition 1.3.9]{Pedersenbook}, as the
main part of the proof of the more general statement saying that
if $0\leq x \leq y$ implies $x^{\beta} \leq y^{\beta}$ for all
$x,y$ in a $C^*$-algebra $A$ and for a positive number $\beta > 1$, then $A$ is commutative. 
In the present terminology, this result says that if the function
$f(t)=t^\beta$ is $A$-monotone on the interval $[0,\infty[$ for
some $\beta > 1$ then the $C^*$-algebra $A$ is commutative. In
1998, W. Wu proved that if $f(t)=e^t$ is $A$-monotone, then the
$C^*$-algebra $A$ is commutative \cite{Wu-art}, by reducing the
proof via involved approximation arguments to the $A$-monotonicity
of the function $t^2$ and then using the Ogasawara's result
\cite{Ogasawara}. In the recent paper by G. Ji and J. Tomiyama
\cite{JiTomiyama} it has been proved that the \mbox{$C^*$-algebra}
is commutative if and only if  all monotone functions are
$A$-monotone, that is $P_1 = P_A$, and also if and only if there
exists a continuous monotone function on the positive axis which
is not matrix monotone of order $2$ but $A$-monotone. If one makes
the use of the operator monotonicity of the $\log t$ function
noted already by C. L{\"o}wner \cite{Loewner-art}, one can deduce
from the $A$-monotonicity assumption for $e^t$, the
$A$-monotonicity of the function $t^{\beta}$ for any 
$\beta > 1$. 
Hence by the above cited result \cite[Proposition 1.3.9]{Pedersenbook}  the $C^*$-algebra $A$ has to be commutative which provides a short proof for 
the above mentioned result of Wu \cite{Wu-art}.

A $C^{*}$-algebra is commutative if and only if all its
irreducible representations are one-dimensional, or in other words
if and only if it is $1$-homogeneous. In this sense both the
$n$-homogeneous $C^*$-algebras and the $n$-subhomogeneous
$C^*$-algebras, that is those $C^*$-algebras having only
$n$-dimensional irreducible representations or respectively only
irreducible representations of dimension less or equal to a positive
integer $n$, are natural generalizations of the class of
commutative $C^*$-algebras.

Using our results on relationship between homogeneity of
$C^*$-algebras and the standard matrix monotonicity scaling of
functions, we obtain an extension of the result of G. Ji and J.
Tomiyama \cite{JiTomiyama} to the $n$-subhomogeneous
$C^*$-algebras.

\begin{theorem} \label{th:intermidiatn-monotcond}
Let $A$ be a $C^*$-algebra. If there exists a pair of positive
integers $(m,n)$ obeying $ n<m$, and such that firstly, every
$A$-monotone function is $n$-monotone, that is $P_A \subset P_n$,
and secondly, there is a function which is at the same time
$A$-monotone, $n$-monotone but not $m$-monotone, that is 
$P_A \cap (P_n \setminus P_m ) \neq \emptyset$, then there exists some
intermediate integer $n\leq j < m$ such that

\begin{itemize}

\item[1)]

every $A$ monotone function is $j$-monotone, that is $P_A = P_j$;

\item[2)]the $C^*$-algebra $A$ is $j$-subhomogeneous.

\end{itemize}

\end{theorem}

If $m=n+1$, then we get the following useful specialization of
Theorem \ref{th:intermidiatn-monotcond}.

\begin{theorem} \label{th:intermidiatnnplus1-monotcond}
Let $A$ be a $C^*$-algebra. If there exists a positive integer
$n$, such that firstly, every $A$-monotone function is
$n$-monotone, that is $P_A \subset P_n$, and secondly, there is a
function which is at the same time $A$-monotone, $n$-monotone but
not $(n+1)$-monotone, that is 
$P_A \cap (P_n \setminus P_{n+1}) \neq
\emptyset$, then $P_A = P_n$ and the $C^*$-algebra $A$ is
$n$-subhomogeneous.

\end{theorem}

\begin{proof}({\em Theorem \ref{th:intermidiatn-monotcond} and
\ref{th:intermidiatnnplus1-monotcond}}) By Lemma
\ref{th:mainlemmareps-omf}, for any $C^*$-algebra $A$ and any
positive integer $k$, only two possibilities are possible, either
$P_A \cap P_k = P_k$ or $P_A \cap P_k = P_A$. If there exists a
positive integer $n$, such that $P_A \subset P_n$, then $A$ is
$n_0$-subhomogeneous and $P_A = P_{n_0}$ by Corollary
\ref{th:maxmonmaxhomog}. We have that 
$P_{n_0} \cap (P_n \setminus P_m ) = P_{A} \cap (P_n \setminus P_m ) 
\neq \emptyset$ and
$P_{n_0} = P_A \subseteq P_n$. Since $P_m \subsetneq P_n$ for all
$m > n$, we get 
$P_m \subseteq P_{n_0} \subseteq P_n$ and $n\leq n_0 < m$ by assumption.

Theorem \label{th:intermidiatnn+1-monotcond} is obtained in the
special case when $m=n+1$. Indeed, in this case 
$P_{n_0}= P_A \subseteq P_n $ and $n \leq n_0 < n+1$. 
Hence $n=n_0$ and $P_A = P_n$. Thus 
$A$ becomes $n$-subhomogeneous.

\end{proof}

\begin{example}
As we have mentioned before  $f(t)= t^\beta \in P_1\setminus P_2$ for
$\beta > 1$.  Hence if $f(t)= t^\beta \in P_A$ for some $\beta > 1$, then 
$P_A \cap (P_1\setminus P_2) \neq \emptyset$, and by Theorem
\label{th:intermidiatnnplus1-monotcond} we get that $P_A = P_1$. 
Hence $A$
is $1$-homogeneous by Theorem \ref{th:monot-homog}, that is all
its irreducible representations are one-dimensional. 
Indeed, this implies that $A$ is commutative. This is 
the essential point of the arguments in 
Ji and Tomiyama \cite{JiTomiyama}, which yields the results of G. K.
Pedersen \cite[Proposition 1.3.9]{Pedersenbook}, T. Ogasawara
\cite{Ogasawara} and W. Wu \cite{Wu-art}.
\end{example}

A complementing assertions to Theorem
\ref{th:intermidiatn-monotcond} is as follows. The corresponding
specialization of Theorem \ref{th:intermidiatn-monotcondempty} for
$m=n+1$ is obtained just by replacing $m$ by $n+1$.

\begin{theorem} \label{th:intermidiatn-monotcondempty}
Let $A$ be a $C^*$-algebra. If for all pairs of positive integers
$(m,n)$ obeying $n<m$ there is no functions that are at the same
time $A$-monotone, $n$-monotone but not $m$-monotone, that is 
$P_A \cap (P_n \setminus P_m ) = \emptyset$, then

\begin{itemize}

\item[1)]

every $A$ monotone function is  operator monotone, that is 
$P_A = P_{\infty}$;

\item[2)] either the set of dimensions of finite-dimensional irreducible representations of
$A$ is unbounded, or $A$ has an infinite-dimensional irreducible
representation.

\end{itemize}

\end{theorem}

\begin{proof}
Suppose that $P_A \neq P_\infty$ in spite of 
$P_A \cap (P_n \setminus P_m ) = \emptyset$. Then by Lemma
\ref{th:mainlemmareps-omf} all irreducible representations of $A$
are finite-dimensional and the set of their dimensions is bounded.
By Theorem \ref{th:monot-homog} there exists a positive integer 
$k \geq n $ such that $P_A = P_{k}$. Since the existence of gaps
asserts that $P_k \subsetneq P_{k+1} \neq \emptyset $, we have 
$P_k \setminus P_{k+1} \neq \emptyset$ and hence 
$P_A \cap (P_k \setminus P_{k+1}) = P_k \cap (P_k \setminus P_{k+1}) = 
(P_k \setminus P_{k+1}) \neq \emptyset $ in contradiction with the
condition of the theorem.

\end{proof}

\begin{remark} The gaps between classes of monotone matrix functions were
addressed in \cite{Donoghuebook}, and more recently, in
\cite{HansenJiTomiyama-art} and \cite{Nayak1}. In \cite{Nayak1}
the "if and only if" extension was obtained of the result on
fractional mapping between classes of matrix monotone functions
from the paper by Wigner and von Neumann \cite{WignervNeumann},
and then it was shown that this extended result yields a proof of
the implication that if $n\ge 2$, then $P_n=P_{n+1}$ implies
$P_n=P_{\infty}$. This can be viewed as a different prove for
specialization of Theorem \ref{th:intermidiatn-monotcondempty} to
the case when $A=B(H)$ and $m = n+1$.
\end{remark}

\begin{remark} In \cite{SparrMathScandart}, a new proof of
L{\"o}wner's theorem on integral
representation of operator monotone functions, different from the
three proofs by L{\"o}wner, Bendat and Sherman, and Karanyi and
Nagi, has been obtained by employing another classes of functions
${\cal M}_n$ in between $P_n$ and $P_{n+1}$. A real-valued
functions $h$ on $(0,\infty)$ is in ${\cal M}_n$ if and only if,
for $a_j \in \mathbb{R}$, $\lambda_j > 0$ and $j = 1,\dots, 2n$, 
the following implication holds:

$$ \left(\sum_{j=1}^{2n} a_j \frac{t \lambda_j - 1}{t+\lambda_j} \geq 0
\mbox{ for } t > 0, \sum_{j=1}^{2n} a_j = 0 \right) \Rightarrow
\left(\quad \sum_{j=1}^{2n} a_j h(\lambda_j) \geq 0\right).
$$ 
As important part of the proof of L{\"o}wner's theorem, it was shown
in \cite{SparrMathScandart} that 
$P_{n+1} \subseteq {\cal M}_n \subseteq P_{n}$ 
for any positive integer $n$. There an explicit
example, showing that $P_2 \setminus {\cal M}_2 \neq \emptyset$,
has been pointed out, thus particularly implying that 
$P_2 \setminus P_3 \neq \emptyset$. 
Proving that 
$P_n \setminus {\cal M}_n \neq \emptyset$ and 
${\cal M}_n \setminus P_{n+1} \neq \emptyset$ for an arbitrary $n$ 
is still an open problem. Motivated by our results, we
feel that the related problem of finding a $C^*$-algebraic
interpretation and perhaps a $C^*$-algebraic generalization of the
spaces ${\cal M}_n$ would be of interest.
\end{remark}

Theorem \ref{th:intermidiatn-monotcond} can be used to obtain the
following unexpected operator monotonicity based characterizations
of subhomogeneous $C^*$-algebras and of dimension for Hilbert
spaces.

Let $g_n (t)= t + \frac{1}{3} t^3 + \dots + \frac{1}{2n-1}
t^{2n-1}$ , where $n$ is some positive integer. 
In \cite{HansenJiTomiyama-art} it was proved that there exists

$\alpha_n > 0$ such that 
$g_n \in P_n([0,\alpha_n[) \setminus P_{n+1}([0,\alpha_n[)$, 
and consequently  
$f_n  = g_n \circ h_{n} \in P_n \setminus P_{n+1}$, 
where $h_n(t)$ is the
M{\"o}bius transformation $h_n(t)= \frac{\alpha_n t}{1+t}$,
operator monotone on $[0,\infty[$, with the inverse 
$h_n^{-1}(t)= \frac{t}{\alpha_n-t}$ operator monotone on 
$[0,\alpha_n[$.

\begin{corollary} \label{cor1:subhomconcrfnc}
If $f_n$ is $A$-monotone function on $[0,\infty [$ for a
$C^*$-algebra $A$, then $A$ is a subhomogeneous $C^*$-algebra, 
such that dimensions of all its
irreducible representations do not exceed $n$.

\end{corollary}

\begin{corollary} \label{cor2:subhomconcrfnc}
If $f_n$ is $B(H)$-monotone for some positive integer $n$ and a
Hilbert space $H$, then $\dim H \leq n$.

\end{corollary}

\begin{proof}(\emph{Corollary \ref{cor1:subhomconcrfnc} and
\ref{cor2:subhomconcrfnc}}) By Remark \ref{rem:PAPkinclud}, 
$P_A \subseteq P_n$ or 
$P_n \subseteq P_A$. If $P_A \subseteq P_n$, then 
$P_A = P_n$ by Theorem  
\ref{th:intermidiatnnplus1-monotcond}. 
Hence $A$ is 
subhomogeneous by Theorem \ref{th:monot-homog}(2).
If $P_A \subseteq P_n$, then there exists 
$k \leq n$ such that $P_A = P_k$, and dimensions of 
irreducible representations do not exceed $n$ 
by Theorem \ref{th:monot-homog}(2). 
In the special case when $A=B(H)$, this property yields $\dim H \leq n$.

\end{proof}

\section{Existence of subalgebras respecting scaling} \label{sec:existsubalg}

In this section we obtain some results on $C^*$-subalgebras and
monotonicity, that can be viewed as a $C^*$-algebraic counterpart
of the standard scaling 
$ M_k \hookrightarrow M_n \hookrightarrow B(H)$, 
$k < n < \dim H = \infty$.

In the the following theorem $CM_m$ means the cone of $M_m$, that is, 
$C_0([0, 1]) \otimes M_m = C_0([0, 1], M_m)$. 
In the proof we will make use of
some results on projective $C^*$-algebras and lifting of relations
in $C^*$-algebras \cite{Loringbook,LoringPedersen1}.

\begin{theorem} \label{th:n-homsubalgexist} Let $A$ be a
$C^*$-algebra.

\begin{itemize}

  \item[1)] If $A$ is a $C^*$-algebra having an $n$-dimensional irreducible
representation for some positive integer $n$, then for any
positive integer $m\leq n$ there exists an $m$-homogeneous
   (presumably nonunital) $C^*$-subalgebra $B$.

  \item[2)] If $A$ has an infinite-dimensional irreducible representation $\pi$, then

   \begin{itemize}

     \item[2a)] For any positive integer $m$ there exists
     a $C^*$-subalgebra $B$ in $A$, such that $B$ is
     $m$-homogeneous.

     \item[2b)] The  $C^*$-algebra has $\infty$-homogeneous  
     $C^*$-subalgebra $B$,
  that is a $C^*$-subalgebra whose all non-zero irreducible representations are
     infinite-dimensional, if and only if $A$ is not residually
     finite-dimensional, that is 
     $$I = \bigcap_{\pi} Ker (\pi) \neq \{0\},$$
     where the intersection is taken over all finite-dimensional
     irreducible representations.

   \end{itemize}

  \item[3)] If the set of dimensions of finite-dimensional irreducible representations of $A$ is unbounded, then for any positive integer $m$ there exists
$m$-homogeneous $C^*$-subalgebra of $A$.

 \end{itemize}

\end{theorem}

\begin{proof}
1) Let $\pi : A \rightarrow B(H)$ be an $n$-dimensional irreducible representation of $A$. Then $\pi(A)$ is isomorphic to $n \times n$ matrix algebra $M_n$.
Let $\{e_{i,j}\}$ be the standard matrix units for $\pi(A)$ 
obtained  from the standard matrix units for 
$n \times n$ matrix algebra via this isomorphism. 
Now for any positive integer 
  $m \leq n$ the elements $a_2 = e_{2,1}, \dots, a_m = e_{m,1}$ satisfy the
  relations of Theorem 10.2.1 in the
Loring's book \cite{Loringbook}, namely, 
$$
   (*) \left\{ \begin{array}{ll}
                 ||a_j|| \leq 1&\\
                 a_ja_k = 0 &(j, k = 2,\dots, m)\\
                 a_j^*a_k = \delta_{j,k} a_2^*a_2 &(j, k = 2, \dots, m)
                 \end{array}
                 \right.
  $$
  Hence by the above cited theorem these elements are lifted to $A$ 
  keeping those relations. Thus there are elements 
  $\bar{a_2}, \dots, \bar{a_m}$ in $A$ satisfying the same relations 
  such that $\pi(\bar{a_j}) = a_j$ for $2 \leq j \leq m$. 

Let $B = C^*(\bar{a_2}, \dots, \bar{a_m})$ be the C*-algebra of $A$
  generated by $\bar{a_2}, \dots, \bar{a_m}$. 

By Proposition 3.3.1 in \cite{Loringbook}, the universal C*-algebra on generators $c_2, \dots, c_m$ satisfying the same relations is isomorphic to $CM_m$ by the map $c_j \mapsto t \otimes e_{j,1}$. Therefore there exists a homomorphism from $CM_m$ onto $B$, and since $CM_m$ is $m$-homogeneous, its image $B$ must be $m$-homogeneous.

2a) Let $\pi$ be an infinite-dimensional irreducible
representation of $A$ on a Hilbert space $H$. Take an
$m$-dimensional projection $p: H \mapsto H$. By Kadison
transitivity theorem, in the form it is stated in Takesaki's book
\cite[Ch.2, Theorem 4.18]{takesaki-bok}, $\pi (A) p = B(H) p$.
Then $pB(H)P = p \pi (A) p \cong M_m$. The restriction of 
$\pi: B \mapsto p B(H) p $ to the $C^*$-subalgebra 
$B=\{a\in A \mid \pi(a) pH \subseteq pH, \pi(a)^* pH \subseteq pH \}$ of $A$
consisting of elements mapped by $\pi$ to operators that together
with their adjoints leave $pH$ invariant, is $n$-dimensional
representation of $B$ on $pH$, and moreover it is irreducible and
surjection, 
since $\pi (B) p = p \pi (B) p = p \pi (A) p = p B(H)p = B (pH)$. 
Repeating the lifting argument from 1) with $k = m$, we
get the $m$-homogeneous $C^*$-subalgebra of $B$ and thus of $A$.

2b)$\Uparrow$: Suppose $A$ is not residually finite-dimensional.
Then
$$I = \bigcap_{
  \begin{array}{c} \pi \in \mbox{ irred.rep.}(A)\\
  \dim \pi < \infty
  \end{array}}
  Ker (\pi) \neq {0}$$
is an ideal and thus is a $C^*$-subalgebra in $A$. Let $\pi$ be a
non-zero irreducible representation of $I$ on a Hilbert space $H$.

Since $I$ is an ideal in $A$ there exists an irreducible
representation $\tilde{\pi}$ of $A$ on the same Hilbert space $H$
extending $\pi$, that is coinciding with $\pi$ on $I$. If $\dim
\tilde{\pi}= \dim \pi = \dim H$ is finite, then $\pi (I) =
\tilde{\pi} (I) = {0}$ by definition of $I$, in contradiction with
assumption that $\pi$ is non-zero. Hence, $\dim \pi = \infty$, and
since $\pi$ has been chosen arbitrary, $B = I$ is an
$\infty$-homogeneous $C^*$-subalgebra of $A$.

2b) $\Downarrow$: Assume that $A$ is residually
finite-dimensional, that is
$$I = \bigcap_{
  \begin{array}{c} \pi \in \mbox{ irred. rep.}(A)\\
  \dim \pi < \infty
  \end{array}}
  Ker (\pi) = {0}.$$

Then any $C^*$-subalgebra $B$ of $A$ has a non-zero
finite-dimensional irreducible representation. Indeed, $A$ being
residually finite-dimensional has sufficiently many
finite-dimensional irreducible representations, that is for any
non-zero $C^*$-subalgebra $B$ there exists an irreducible
finite-dimensional representation $\tilde{\pi}$ of $A$ such that
$\tilde{\pi} (B) \neq \{0\}$. 
Since $\tilde{\pi} (A) \cong M_{\dim \tilde{\pi}}$, 
it holds that 
$\tilde{\pi}(B) = M_{k_1} \oplus \dots \oplus M_{k_l}$ 
is the direct sum of full matrix algebras.

Then cutting down onto one of the summands by a central projection
$p_j$ yields a non-zero finite-dimensional irreducible
representation $\pi: B \rightarrow p_j \tilde{\pi} (B) p_j$ of
$B$. Hence, there is no $\infty$-homogeneous $C^*$-subalgebras in
$A$, if $A$ is residually finite-dimensional.

3) If the set of dimensions of finite-dimensional irreducible
representations of $A$ is unbounded, then for any positive integer
$m$ there exists an irreducible representation of dimension 
$n > m$. As in 1), the elements $a_2 = e_{21}, \dots, a_m = e_{m1}$
of $M_{n}$ can be lifted to the elements of $A$, which generate an
$m$-homogeneous $C^*$-subalgebra in $A$.

\end{proof}

Combining Lemma \ref{th:mainlemmareps-omf} with Theorem
\ref{th:n-homsubalgexist} we obtain the following result.

\begin{theorem} \label{th:subalgexistn-monoton} Let $A$ be a
$C^*$-algebra.

\begin{itemize}

   \item[1)] If $A$ is a $C^*$-algebra having $n$-dimensional irreducible
representation for some positive integer $n$, then for any
positive integer $m\leq n$ there exists a $C^*$-subalgebra $B$,
such that $P_B = P_m$.

   \item[2)] If $A$ has an infinite-dimensional irreducible representation $\pi$, or
    the set of dimensions of finite-dimensional irreducible representations
    of $A$ is unbounded, then for any positive integer $m$ there exists
    $C^*$-subalgebra $B$ of $A$ such that $P_B = P_m$.

   \item[3)] If $A$ is not residually finite dimensional, then there exists a $C^*$-subalgebra
   $B$ such that $P_B = P_\infty $.

 \end{itemize}

\end{theorem}

\begin{example}
Let $A = C(X) \rtimes_\sigma \mathbb{Z}$ be the transformation
group (crossed product) $C^*$-algebra associated to a dynamical
system $\Sigma = (X, \sigma)$ consisting of a homeomorphism
$\sigma$ on a compact Hausdorff metric space $X$, and let
$Per(\Sigma)$ denote the set of all periodic points of $\Sigma$,
and $\overline{Per(\Sigma)}$ be its closure in $X$. In
\cite[Theorem 4.6]{TomSeoulLN1} it was shown that $A$ is
residually finite-dimensional if and only if
$\overline{Per(\Sigma)} = X$. Hence, by Theorem
\ref{th:n-homsubalgexist} and Theorem
\ref{th:subalgexistn-monoton}, if $\overline{Per(\Sigma)} \neq X$,
then $A = C(X) \rtimes_\sigma \mathbb{Z}$ contains some
$C^*$-subalgebra $B$ such that all irreducible representations of
$B$ are infinite-dimensional and $P_B = P_\infty $. The group
$C^*$-algebra 
$C^*({\cal H}) = C(\mathbb{T}^2) \rtimes_\sigma\mathbb{Z}$ 
of the three-dimensional Heisenberg group mensioned in
Example \ref{ex:Heisgroup} has property 
$\overline{Per(\Sigma)} = X$, and is therefore a residually finite-dimensional
$C^*$-algebra. By Theorem \ref{th:n-homsubalgexist}, any
$C^*$-subalgebra of $C^*({\cal H})$ possesses finite-dimensional
irreducible representation. Moreover $C^*({\cal H})$ has
infinite-dimensional irreducible representations (for example
those induced by aperiodic points), and hence by Theorem
\ref{th:n-homsubalgexist} and Theorem
\ref{th:subalgexistn-monoton}, for any positive integer $m$ it has
to have $m$-homogeneous $C^*$-subalgebra $B$ such that 
$P_B = P_m$. It could be interesting to construct explicitly such
subalgebras in $C^*({\cal H})$.

\end{example}

\begin{example}
The $C^*$-algebra of compact operators $K(H)$ on an
infinite-dimensional Hilbert space $H$ is not residually
finite-dimensional. All irreducible representations of $K(H)$ are
infinite-dimensional except for the zero representation. 
Hence 
$$I
= \bigcap_{
  \begin{array}{c} \in \mbox{ irred.rep.}(A)\\
  \dim \pi < \infty
  \end{array}}
  Ker (\pi) = K(H), $$

and $K(H)$ itself can be taken as an example of a
$\infty$-homogeneous $C^*$-subalgebra of $K(H)$ such that

$P_{K(H)} = P_\infty$.

\end{example}

\begin{example}
Let $A$ be a simple $C^*$-algebra, that is a $C^*$-algebra with no
non-zero closed ideals, and assume that $A$ has
infinite-dimensional irreducible representation, thus implying
that $P_A = P_\infty$ by 4) of Lemma \ref{th:mainlemmareps-omf}.
In this case every non-zero 
irreducible representation of A is naturally infinite-dimensional. 
In particular $A$ is not residually finite-dimensional, and so has
at least one $\infty$-homogeneous  $C^*$-subalgebra $B$ such that
$P_B = P_\infty$, namely $B = A$. 
The $C^*$-algebra of compact operators $K(H)$ on
an infinite-dimensional Hilbert space $H$, the Cuntz
$C^*$-algebras ${\cal O}_n$ and the irrational rotation
$C^*$-algebra are simple $C^*$-algebras which have
infinite-dimensional irreducible representations. The question
arising from these observations is whether any infinite-dimensional 
simple $C^*$-algebra contains a proper (different from the whole
algebra) $\infty$-homogeneous $C^*$-subalgebra, and how to find
and classify such $C^*$-subalgebras for the specific examples
where such subalgebras exist.

\end{example}

The following question is suggested by Theorems
\ref{th:n-homsubalgexist}, \ref{th:subalgexistn-monoton} and
\ref{th:intermidiatn-monotcond}:  if $A$ is $n$-matrix monotone,
that is $P_A = P_n$, and its $C^*$-subalgebra $B$ is $k$-matrix
monotone for some $k < n$, that is $P_B = P_k$, then is it true
that for any $l$ between $n$ and $k$ there exists a
$C^*$-subalgebra $C$ containing $B$ for which $C$ is $l$-monotone ?

This question is closely related to the following question
concerning representations of $C^*$-algebras: For $A$ having
irreducible representation of finite dimension $n$, its subalgebra
$B$, sub-homogeneous of degree $k < n$, and any integer $l$
between $k$ and $n$, can we find a sub-homogeneous subalgebra of
degree $l$ including $B$ ?

As an easy counter example to both questions 
let $A = M_4$ and $B=M_2 \oplus M_2$ 
be the direct sum of two $M_2$, which is isomorphic to
a maximal $C^*$-subalgebra of $M_4$ obtained by placing the direct
summands as diagonal blocks. In this case $P_A = P_4$ and
$P_B=P_2$, but there is no $C^*$-subalgebra $C$ of $A$ containing
$B$ with $P_C = P_3$. At the same time, if 
$B=M_2 \oplus M_1 \oplus M_1$ imbedded in $A=M_4$ as diagonal blocks with
non-increasing order of dimensions, then $P_A = P_4$, $P_B=P_2$
and there exists a $C^*$-subalgebra $C$ of $A$ containing $B$ with
$P_C = P_3$. Namely, one can take $C=M_3$ imbedded in $A=M_4$ by
placing the direct summand as the diagonal block containing the
$M_2 \oplus M_1$ part of $B$, and putting $M_1$ in the remaining
diagonal spot. Thus even in the case of the matrix algebra $M_n$
we have to consider the location of its $C^*$-subalgebra $B$ in
$A$ to be able to assert the existence of the intermediate
$C^*$-subalgebra $C$. In other words one needs many experiments in
concrete $C^*$-algebras to clarify how the gaps appear depending
on the nature of the imbeddings. We feel that the monograph
\cite{SCPowers-book}, on subalgebras of $C^*$-algebras and the
limit algebras of inclusion sequences, discussing the importance
of the nature of the inclusions, contains results which could be
of interest in this respect.

As we discussed before for a $C^*$-subalgebra $B$ of a
$C^*$-algebra $A$, the equality $P_B=P_1$ holds if and only if $B$
is commutative (abelian). Maximal abelian $C^*$-subalgebras are
important for understanding representations and structure of a
$C^*$-algebra. Closely related to the previous discussion is the
following problem. Let $B$ be a maximal abelian $C^*$-subalgebra
of a $C^*$-algebra $A$. What are the "allowed" positive integers
$j$ for which there exists a $C^*$-subalgebra $C$ of $A$
containing $B$ such that $P_C = P_j$ ?

Beyond the class of full matrix algebras, there are many important
examples of infinite decreasing inclusion sequences of
$C^*$-subalgebras $A_1 \hookleftarrow A_2 \hookleftarrow A_3
\hookleftarrow \dots $. If for some positive integer $k$ there
exists a positive integer $n$ such that $P_{A_k} = P_n$, then for
all positive integers $j \geq k$ there exists a positive integer
$l_j \leq n$ such that $P_{A_j} = P_{l_j}$, and moreover 
$l_{j+1} \leq l_j $ for all $j\geq n$. 
So the inclusion non-increasing
sequence of function spaces $P_{A_k}$ stabilizes at some positive
integer $s\leq n$, which means that 
$P_{A_j} =  P_s$ for all $j \geq s$.  
The first question arising from these considerations is
whether it is possible to have $ P_{A_{\infty}} = P_t$ for some
$t < s$, where 
$A_{\infty} = \cap_{r \in {\mathbb{N}\setminus
\{0\}}} A_r$. The second question is concerned with $P_{\infty}$.
Suppose that $P_{A_k} = P_\infty$ for all positive integers $k$.
Is it possible to have a decreasing sequence satisfying
$P_{A_{\infty}} = P_\infty$, and if it is possible what are then
the properties of sequences of $C^*$-subalgebras and their
imbeddings leading to this situation ?

%{\small \it Acknowledgements. }


\begin{thebibliography}{99}

\bibitem{Donoghuebook} W. F. Donoghue, Jr., 
\emph{Monotone matrix functions an
analytic continuation}, Springer-Verlag, New York, 1974.

\bibitem{EffrosRuanbook} E. Effros and Z-J. Ruan, \emph{Operator Spaces},
London Math.Soc. Monographs 23, 2000, Oxford Sci.Pub.

\bibitem{Fellart1} J. M. G. Fell, \emph{The structure of algebras of
operator fields}, Acta Math. \textbf{106}, (1961), 233-280.

\bibitem{HansenJiTomiyama-art} F. Hansen, G. Ji, J. Tomiyama, \emph{Gaps
between classes of matrix monotone functions}, to appear in Bul.
London Math. Soc.

\bibitem{HiaiYanagibook} F. Hiai, K.Yanagi, \emph{Hilbert spaces and linear operators},
Makino Pub.Ltd., 1995.

\bibitem{JiTomiyama} G. Ji, J. Tomiyama, \emph{On characterizations of commutativity
 of $C^*$-algebras}, Proc. Amer. Math. Soc.,
 posted on March 25, 2003, S 0002-9939(03)06991-0
(to appear in print).

\bibitem{LoringPedersen1} T. A. Loring, G. K. Pedersen,
\emph{ Projectivity, transitivity and AF-telescopes}, Trans. Amer.
Math. Soc. \textbf{ 350}, 11 (1998), 4313--4339.

\bibitem{Loringbook} T. A. Loring, \emph{Lifting solutions to perturbing
problems in $C^*$-algebras}, Fields Institute Monographs, AMS,
1997.

\bibitem{Loewner-art} C. L{\"o}wner, \emph{ \"{U}ber monotone
Matrixfunktionen}, Math. Z., \textbf{ 38} (1934) 177-216.
\bibitem{Nayak1} S. Nayak,
\emph{Monotone matrix functions of successive orders},
 Proc. Amer. Math. Soc.,  posted on July 17, 2003, PII S 0002-9939(03)07218-6
(to appear in print).

\bibitem{Ogasawara} T. Ogasawara, \emph{A theorem on operator algebras},
J. Sci. Hiroshima Univ., \textbf{ 18} (1955), pp 307-309.

\bibitem{Pedersenbook} G. K. Pedersen, $C^*$-algebras and their automorphism Groups,
Academic press, 1979.

\bibitem{SCPowers-book} S. C. Power, \emph{Limit algebras: an
introduction to subalgebras of $C^*$-algebras}, Pitman Research
Notes in Mathematics Series, 278,  Longman Scientific \&
Technical, Harlow, 1992.

\bibitem{SparrMathScandart} G. Sparr,
\emph{A new proof of L\"owner's theorem on monotone matrix
functions}, Math. Scand. \textbf{47}, 2 (1980), 266-274.

\bibitem{SilTomTypeIart}
S. D. Silvestrov, J. Tomiyama, Topological dynamical systems of
type I, Expo. Math. \textbf{20}, 2 (2002), 117--142.

\bibitem{Stinespring-art} W. F. Stinespring,  \emph{Positive functions on $C^*$-algebras},
 Proc. Amer. Math. Soc. \textbf{ 6} (1955),  211-216.

\bibitem{takesaki-bok}  M. Takesaki, \emph{ Theory of operator algebras I},
Springer-Verlag New York Inc., 1979. \\

\bibitem{TomTohokuJ62} J. Tomiyama, \emph{Topological representations of $C^*$-algebras}, Tohoku Math.J., 
\textbf{14} (1962), 187-204. 

\bibitem{TomSeoulLN1} J. Tomiyama,
\emph{The Interplay between Topological Dynamics and Theory of
$C^*$-algebras}, Lecture Note No.2, Global Anal. Research Center,
Seoul, 1992.

\bibitem{TomTak61} J. Tomiyama, M.Takesaki, 
\emph{Applications of fibre bundles to the certain 
class of C*-algebras}, Tohoku Math. J., \textbf{13} (1961),498-523

\bibitem{WignervNeumann}
E. P. Wigner, J. v. Neumann, \emph{ Significance of Loewner's
theorem in the quantum theory of collisions}, Ann. of Math. (2),
\textbf{59} (1954), 418--433.

\bibitem{Wu-art} W. Wu,  \emph{An order characterization of
commutativity for $C^{\ast}$-algebras},
 Proc. Amer. Math. Soc., \textbf{129}(2001),  983--987.

\end{thebibliography}
\end{document}